\documentclass[12pt]{amsart}

\usepackage{amsmath}
\usepackage{amsfonts}
\usepackage{amssymb}
\usepackage{graphicx}
\usepackage{hyperref}
\usepackage{mathrsfs}
\usepackage{pb-diagram}
\usepackage{epstopdf}
\usepackage{amsmath,amsfonts,amsthm,enumerate,amscd,latexsym,curves}
\usepackage{bbm}
\usepackage[mathscr]{eucal}
\usepackage{epsfig,epsf,xypic,epic}
\usepackage{caption}
\usepackage{subcaption}

\newtheorem{thm}{Theorem}[section]

\newtheorem{conj}[thm]{Conjecture}

\numberwithin{equation}{section}

\def\bZ{\mathbb{Z}}

\def\bR{\mathbb{R}}
\def\bC{\mathbb{C}}
\def\bP{\mathbb{P}}

\begin{document}

\title[SYZ and its impact]{The Strominger-Yau-Zaslow conjecture and its impact}
\author[K. Chan]{Kwokwai Chan}
\address{Department of Mathematics\\ The Chinese University of Hong Kong\\ Shatin\\ Hong Kong}
\email{kwchan@math.cuhk.edu.hk}

\begin{abstract}
This article surveys the development of the SYZ conjecture since it was proposed by Strominger, Yau and Zaslow in their famous 1996 paper \cite{SYZ96}, and discusses how it has been leading us to a thorough understanding of the geometry underlying mirror symmetry.
\end{abstract}

\maketitle


\section{Mirror symmetry before SYZ}

When Yau proved the famous Calabi conjecture \cite{Yau77, Yau78} in 1976, most mathematicians, probably Yau himself included, would not have imagined that {\em Calabi-Yau manifolds} were going to play such a pivotal and indispensable role \cite{CHSW85, Strominger86} in {\em string theory} - a modern physical theory which is aimed at unifying general relativity and quantum field theory. Nor would people speculate that the resulting interaction between geometry and physics would eventually lead to the important discovery of {\em mirror symmetry} - a subject that has been drastically influencing the development of many branches of mathematics for more than two decades.

The story began in the late 1980's when Dixon \cite{Dixon88} and Lerche, Vafa and Warner \cite{LVW89} suggested that theoretically two different Calabi-Yau manifolds could give rise to identical physics. This surprising prediction was soon verified by Greene and Plesser \cite{Greene-Plesser90} and Candelas, Lynker and Schimmrigk \cite{CLS90} when they independently constructed pairs of Calabi-Yau manifolds which exhibit an interchange of Hodge numbers. We call these {\em mirror pairs}. In their remarkable 1991 paper \cite{COGP91}, Candelas, de la Ossa, Green and Parkes carried out an even more astonishing calculation which led to a prediction about the numbers of rational curves on the quintic 3-fold in $\bP^4$. Mathematicians were particularly intrigued by their prediction because it went way beyond what algebraic geometers could achieve at that time.

This has triggered the development of many important subjects such as Gromov-Witten theory, and finally culminated in the proofs of {\em mirror theorems} by Givental \cite{Givental96, Givental98} and Lian-Liu-Yau \cite{LLY-I, LLY-II, LLY-III, LLY-IV} independently, which in particular verified the predictions by Candelas et al. for the quintic 3-fold. This is certainly a magnificent achievement. However, all the proofs rely on the geometry of the ambient toric varieties which contain the Calabi-Yau manifolds, and in particular, they do not provide an intrinsic way to understand the geometry of mirror symmetry.

The first intrinsic mathematical formulation of mirror symmetry was Kontsevich's {\em Homological Mirror Symmetry} (HMS) conjecture, proposed in his 1994 ICM address \cite{Kontsevich-ICM94}. In string theory, a Calabi-Yau manifold $X$ determines two topological string theories: the {\em A-model} and {\em B-model}, which depends on the symplectic and complex geometry of $X$ respectively \cite{Vafa92, Witten92}. From this perspective, mirror symmetry predicts that if $X$ and $\check{X}$ are a mirror pair of Calabi-Yau manifolds, then there is an isomorphism between the A-model of $X$ and the B-model of $\check{X}$, and vice versa. The above enumerative predictions about the quintic 3-fold is one of the many interesting manifestations of this bigger picture.

Kontsevich's HMS conjecture formulates mirror symmetry succinctly as an equivalence between the Fukaya category of Lagrangian submanifolds in $X$ (A-model) and the derived category of coherent sheaves on the mirror $\check{X}$ (B-model). His conjecture is both deep and elegant, and is expected to imply the enumerative predictions by mirror symmetry. Nevertheless, it does not indicate how such and equivalence can be found, nor does it tell us how to construct the mirror of a given Calabi-Yau manifold.

\section{Formulation of the SYZ conjecture}

In the summer of 1996, Strominger, Yau and Zaslow \cite{SYZ96} made a ground-breaking proposal which gave the first geometric explanation for mirror symmetry:
\begin{conj}[The SYZ conjecture \cite{SYZ96}]\label{conj:SYZ}
Suppose that $X$ and $\check{X}$ are Calabi-Yau manifolds mirror to each other. Then
\begin{itemize}
\item[(i)]
both $X$ and $\check{X}$ admit special Lagrangian torus fibrations with sections $\mu:X\to B$ and $\check{\mu}:\check{X}\to B$ over the same base:
\begin{equation*}
\xymatrix{
X \ar[d]_{\mu} & {\check{X}} \ar[d]^{\check{\mu}}\\
{B} \ar@{=}[r] & {B}
}
\end{equation*}

\item[(ii)]
the fibrations $\mu:X\to B$ and $\check{\mu}:\check{X}\to B$ are fiberwise dual to each other in the sense that if the fibers $\mu^{-1}(b) \subset X$ and $\check{\mu}^{-1}(b) \subset \check{X}$ over $b \in B$ are nonsingular, then they are dual tori; and

\item[(iii)]
there exist fiberwise Fourier(-Mukai)--type transforms which are responsible for the interchange between the symplectic-geometric (resp. complex-geometric) data on $X$ and the complex-geometric (resp. symplectic-geometric) data on $\check{X}$.
\end{itemize}
\end{conj}

In a nutshell, this is saying that the mysterious mirror phenomenon is simply a Fourier transform! This remarkable and far-reaching conjecture not only provides a beautiful geometric explanation to mirror symmetry, but also suggests that a mirror partner of any given Calabi-Yau manifold $X$ can be constructed by fiberwise dualizing a special Lagrangian torus fibration on $X$ (or {\em $T$-duality}). It immediately attracted much attention from both mathematicians and physicists, and has lead to a flourishing of research work aiming at either solving the conjecture or applying it to understand the geometry underlying mirror symmetry.

Before going on, let us go through briefly the heuristic arguments behind the SYZ conjecture. First of all, a key feature in string theory is the existence of {\em Dirichlet branes}, or D-branes. Physical arguments suggest that D-branes in the B-model (or simply {\em B-branes}) are coherent sheaves over complex subvarieties while D-branes in the A-model (or {\em A-branes}) are special Lagrangian submanifolds equipped with flat connections. As mirror symmetry predicts an isomorphism between the A-model of $X$ and the B-model of $\check{X}$, the moduli space of an A-brane on $X$ should be identified with the moduli space of the mirror B-brane on $\check{X}$.

Now, points on $\check{X}$ can certainly be regarded as B-branes. And as $\check{X}$ itself is the moduli space these B-branes, it should be identified with the moduli space of certain A-branes $(L,\nabla)$ on $X$, where $L\subset X$ is a special Lagrangian submanifold and $\nabla$ is a flat $U(1)$-connection on $L$. Also, since $\check{X}$ is swept by its points, $X$ should be swept by these special Lagrangian submanifolds $L$ as well. By McLean's theorem \cite{McLean98}, the moduli space of a special Lagrangian submanifold $L \subset X$ is unobstructed and modelled on $H^1(L;\bR)$, while the moduli space of flat $U(1)$-connections (modulo gauge) on $L$ is given by $H^1(L;\bR)/H^1(L;\bZ)$. Therefore, in order to match the dimensions, we should have $\textrm{dim }H^1(L;\bR) = \textrm{dim}_\bC \check{X} = n$. Hence $X$ should admit a special Lagrangian torus fibration
$$\mu:X \to B.$$
Moreover, the manifold $\check{X}$ itself can be regarded as a B-brane whose moduli space is a singleton and it intersects each point in $\check{X}$ once, so the corresponding A-brane should give a special Lagrangian section $\sigma$ to $\mu$ with $H^1(\sigma;\bR) = 0$. In particular, the base $B$ should have first Betti number $b_1=0$.

Applying the same argument to $\check{X}$ yields a special Lagrangian torus fibration with section
$$\check{\mu}:\check{X} \to \check{B}.$$
Now for a torus fiber $L_b := \mu^{-1}(b) \subset X$, its dual $L_b^\vee$ can be viewed as the moduli space of flat $U(1)$-connections on $L$ which, under mirror symmetry, correspond to points in $\check{X}$. This shows that $L_b^\vee$ is a submanifold in $\check{X}$. With more elaborated arguments, one can see that $L_b^\vee$ can in fact be identified with a special Lagrangian torus fiber of $\check{\mu}$, and hence deduce that $\mu$ and $\check{\mu}$ are fibrations over the same base which are fiberwise dual to each other.

Notice that we have a transform carrying special Lagrangian torus fibers in $X$ (A-branes) to points in $\check{X}$ (B-branes). This is an instance of a fiberwise Fourier(-Mukai)--type transform. More generally, there should exist geometric Fourier transforms mapping symplectic-geometric data on $X$ to complex-geometric data on $\check{X}$. We call these {\em SYZ transforms}. In the original SYZ paper \cite{SYZ96}, it was inferred that the behavior of the Ricci-flat metrics on the mirror $\check{X}$ should differ from the {\em semi-flat Calabi-Yau metrics}, constructed earlier by Greene, Shapere, Vafa and Yau in an important paper \cite{GSVY90}, by contributions from {\em instanton corrections}. As we shall see, a key step in the investigation of mirror symmetry is to understand these corrections, which should come from higher Fourier modes of the SYZ transforms.

\section{Semi-flat SYZ}

In case the special Lagrangian torus fibrations do not admit any singular fibers, the SYZ picture is particularly nice. McLean's classic results \cite{McLean98} give us two naturally defined {\em integral affine structures}\footnote{An integral affine structure on a manifold is an atlas of charts whose transition maps are all integral affine linear transformations.} on the base manifold $B$ of a special Lagrangian torus fibration: the {\em symplectic} and {\em complex} affine structures, and mirror symmetry can be explained neatly via these structures. More specifically, a normal vector field $v$ to a fiber $L_b := \mu^{-1}(b)$ determines a 1-form $\alpha := -\iota_v \omega \in \Omega^1(L_b;\bR)$ and an $(n-1)$-form $\beta := \iota_v \textrm{Im }\Omega \in \Omega^{n-1}(L_b;\bR)$, where $\omega$ and $\Omega$ are the K\"ahler form and holomorphic volume form on $X$ respectively. McLean \cite{McLean98} proved that the corresponding deformation is special Lagrangian if and only if both $\alpha$ and $\beta$ are closed. By identifying $TB$ with $H^1(L_b;\bR)$ using the cohomology class of $\alpha$, we get the symplectic affine structure on $B$, while identifying $TB$ with $H^{n-1}(L_b;\bR)$ using the cohomology class of $\beta$ gives us the complex affine structure on $B$. We also have the {\em McLean metric} defined by
$$g(v_1,v_2) := -\int_{L_b} \iota_{v_1}\omega \wedge \iota_{v_2}\textrm{Im }\Omega.$$

In his illuminating paper \cite{Hitchin97}, Hitchin explains how these structures are all related through the {\em Legendre transform}. If we denote by $x_1,\ldots,x_n$ the local affine coordinates on $B$ with respect to the symplectic affine structure, then locally the McLean metric can be written as the Hessian of a convex function $\phi$ on $B$, i.e.
$g\left(\frac{\partial}{\partial x_i}, \frac{\partial}{\partial x_j}\right) = \frac{\partial^2\phi}{\partial x_i\partial x_j}$.
Furthermore, setting $\check{x}_i := \partial\phi/\partial x_i$ ($i=1,\ldots,n$) gives precisely the local affine coordinates on $B$ with respect to the complex affine structure, and if
$$\check{\phi} := \sum_{i=1}^n \check{x}_i x_i - \phi(x_1,\ldots,x_n)$$
is the Legendre transform of $\phi$, then we have $x_i = \partial\check{\phi}/\partial\check{x}_i$ and
$g\left(\frac{\partial}{\partial \check{x}_i}, \frac{\partial}{\partial \check{x}_j}\right) = \frac{\partial^2\check{\phi}}{\partial \check{x}_i\partial \check{x}_j}$.

If additionally we assume that the fibration $\mu:X \to B$ admits a Lagrangian section, then a theorem of Duistermaat \cite{Duistermaat80} implies that there are global {\em action-angle coordinates} so that we can write
$$X = T^*B/\Lambda^\vee,$$
where the lattice $\Lambda^\vee \subset T^*B$ is locally generated by $dx_1,\ldots,dx_n$, and $\omega$ can be identified with the canonical symplectic form
$$\omega = \sum_{i=1}^n dx_i \wedge du_i$$
on $T^*B/\Lambda^\vee$. Here $u_1,\ldots,u_n$ are the fiber coordinates on $T^*B$.

In this case, the mirror of $X$ is simply given by
$$\check{X} := TB/\Lambda,$$
where the lattice $\Lambda \subset TB$ is locally generated by $\partial/\partial x_1,\ldots,\partial/\partial x_n$. The quotient $\check{X}$ has a natural complex structure whose holomorphic coordinates are given by $z_i := \exp(x_i + \mathbf{i}y_i)$, where $y_1,\ldots,y_n$ are the fiber coordinates on $TB$ dual to $u_1,\ldots,u_n$. This constructs the mirror of $X$ as a complex manifold with a nowhere vanishing holomorphic volume form
$$\check{\Omega} := d\log z_1 \wedge \cdots \wedge d\log z_n.$$
Moreover, there is an explicit fiberwise Fourier--type transform, which we call the {\em semi-flat SYZ transform} $\mathcal{F}^\textrm{semi-flat}$, that carries $\exp \mathbf{i}\omega$ to $\check{\Omega}$; see \cite[Section 2]{Chan-Leung10b} for more details.

Now if we switch to the complex affine structure on $B$, then we get a symplectic structure on $\check{X}$ which is compatible with its complex structure so that the mirror $\check{X}$ becomes a K\"ahler manifold. Furthermore, if the function $\phi$ above satisfies the real Monge-Amp\`ere equation
$$\textrm{det}\left(\frac{\partial^2\phi}{\partial x_i\partial x_j}\right) = \textrm{constant},$$
then we obtain $T^n$-invariant Ricci-flat metrics on both $X$ and its mirror $\check{X}$. The induced metric on $B$ is called a {\em Monge-Amp\`ere metric} and $B$ is called a {\em Monge-Amp\`ere manifold}. This links mirror symmetry to the study of real Monge-Amp\`ere equations and affine K\"ahler geometry, where Cheng and Yau had made substantial contributions \cite{Cheng-Yau77, Cheng-Yau82, Cheng-Yau86} before even mirror symmetry was discovered. The construction of Monge-Amp\`ere metrics on affine manifolds with singularities has since been an important question in both affine geometry and the study of the SYZ conjecture. The highly nontrivial works of Loftin, Yau and Zaslow \cite{Loftin05, LYZ05} constructed such metrics near the ``Y'' vertex, a typical type of singularity in the 3-dimensional case. But other than this, not much is known.

So the SYZ conjecture indeed paints an appealing picture for mirror symmetry in the semi-flat case; many more details on semi-flat SYZ mirror symmetry were worked out by Leung in \cite{Leung05}. Unfortunately, this nice picture can hold true only at the large complex structure/volume limits where all instanton corrections are suppressed. Away from the limits, special Lagrangian fibrations will have singular fibers and the mirror can no longer be obtained simply by dualizing a fibration.


\section{Constructing SYZ fibrations}


Right after the introduction of the SYZ conjecture in 1996, a great deal of effort was input into constructing special Lagrangian torus fibrations, or {\em SYZ fibrations}, on Calabi-Yau manifolds.\footnote{In \cite{Gross98, Gross99}, instead of constructing such fibrations, Gross assumed their existence and deduced interesting consequences which were predicted by mirror symmetry.} Zharkov \cite{Zharkov00} was the first to construct {\em topological} torus fibrations on Calabi-Yau hypersurfaces in a smooth projective toric variety $\bP_\Delta$. This includes the important example of the quintic 3-fold. Zharkov obtained his fibrations by deforming the restriction of the moment map on $\bP_\Delta$ to the boundary $\partial \Delta$ of the moment polytope to a nearby smooth Calabi-Yau hypersurface.

Using similar ideas and introducing a gradient-Hamiltonian flow, W.-D. Ruan constructed {\em Lagrangian} torus fibrations on quintic 3-folds in a series of papers \cite{Ruan01, Ruan02, Ruan03}. He also carried out an important computation of the monodromy of the fibrations, which was later used by Gross \cite{Gross-Inventiones} to work out a topological version of SYZ mirror symmetry for Calabi-Yau manifolds including the quintic 3-fold. There is also a related work of Mikhalkin \cite{Mikhalkin04} which produces smooth torus fibrations on hypersurfaces in toric varieties by applying tools from tropical geometry.

In general, the construction of special Lagrangian fibrations, or even special Lagrangian submanifolds, is a very difficult problem. One promising approach in constructing special Lagrangians is using the {\em mean curvature flow}. Thomas \cite{Thomas01} formulated a notion of stability for classes of Lagrangian submanifolds with Maslov index zero in a Calabi-Yau manifold, which should be mirror to the stability of holomorphic vector bundles. In \cite{Thomas-Yau02}, Thomas and Yau conjectured that there should exist a unique special Lagrangian representative in a Hamiltonian isotopy class if and only if the class is stable, and that such a representative could be obtained by the mean curvature flow where long time existence should hold. They further proposed a Jordan-H\"older--type decomposition for special Lagrangian submanifolds and related this to formation of singularities in the mean curvature flow. Their proposals and conjectures have a big influence on the development of Calabi-Yau geometry and the SYZ conjecture. There has been a lot of advances in this area \cite{Chen-Li01, Chen-Li04, JLT10, Lee-Wang09, Neves07, Neves10, Neves13, Smoczyk02, Smoczyk04, Smoczyk-Wang02, Wang01, Wang02, Wang03, Wang08a}; see the excellent survey articles \cite{Wang08} and \cite{Neves11} and references therein for more details. Unfortunately, at the time of writing, it is still unknown whether there exists a special Lagrangian torus fibration on the quintic 3-fold.


Noncompact examples of special Lagrangian fibrations are much easier to come by. Harvey and Lawson's famous paper on calibrated geometries \cite{Harvey-Lawson82} gave the simplest of such examples: the map defined by
\begin{align*}
f: \bC^3 & \to \bR^3,\\
(z_1,z_2,z_3) & \mapsto \left(\textrm{Im}(z_1z_2z_3),|z_1|^2-|z_2|^2,|z_1|^2-|z_3|^2\right),
\end{align*}
is a special Lagrangian fibrations whose fibers are invariant under the diagonal $T^2$-action on $\bC^3$. This example was later largely generalized by independent works of Goldstein \cite{Goldstein01} and Gross \cite{Gross01}. They constructed explicit special Lagrangian torus fibrations on any {\em toric} Calabi-Yau $n$-fold (which are necessarily noncompact) via the compatible $T^{n-1}$-action which preserves the natural holomorphic volume form. The discriminant loci of these examples are of real codimension two and can be described explicitly.

Another set of noncompact examples, which has a historic impact on the development of SYZ mirror symmetry and special Lagrangian geometry, was discovered by Joyce \cite{Joyce03}. It was once believed that special Lagrangian fibrations would always be smooth and hence have codimension two discriminant loci. But the examples of Joyce showed that this is unlikely the case. He constructed explicit $S^1$-invariant special Lagrangian fibrations which are only {\em piecewise smooth} and have {\em real codimension one discriminant loci}. The set of singular points of such a fibration is a Riemann surface and its amoeba-shaped image gives the codimension one discriminant locus. Joyce also argued that his examples exhibited the generic behavior of discriminant loci of special Lagrangian fibrations.


This pioneering work of Joyce significantly deepens our understanding of the singularities of special Lagrangian fibrations, and at the same time forces us to rethink about the formulation of the SYZ conjecture. Originally, we expect that a mirror pair of Calabi-Yau manifolds $X$ and $\check{X}$ should have special Lagrangian torus fibrations to the same base $B$ so that their discriminant loci coincide. But the examples of Joyce demonstrate that while the discriminant locus on one side may be of codimension one, that on the other side can be of codimension two. The best that one can hope for is that both $X$ and $\check{X}$ admit special Lagrangian torus fibrations to the same base $B$ and as one approaches the large complex structure limits on both sides, the discriminant loci of these fibrations converge to the same limit which is of codimension two.

More precisely, let $\mathfrak{X} \to D$ and $\check{\mathfrak{X}} \to D$ be maximally unipotent degenerations of Calabi-Yau manifolds mirror to each other, where $D$ is the unit disk and $0\in D$ corresponds to large complex structure limits for the mirror pair. We choose a sequence $\{t_i\} \subset D$ converging to 0, and let $g_i$ and $\check{g}_i$ be Ricci-flat metrics on $\mathfrak{X}_{t_i}$ and $\check{\mathfrak{X}}_{t_i}$ respectively normalized so that they have fixed diameters $C$. Then we expect that
\begin{itemize}
\item[(i)]
there are convergent subsequences of $(\mathfrak{X}_{t_i}, g_i)$ and $(\check{\mathfrak{X}}_{t_i}, \check{g}_i)$ converging (in the Gromov-Hausdorff sense) to metric spaces $(B_\infty, d_\infty)$ and $(\check{B}_\infty, \check{d}_\infty)$ respectively;
\item[(ii)]
the spaces $B_\infty$ and $\check{B}_\infty$ are affine manifolds with singularities which are both homeomorphic to $S^n$;
\item[(iii)]
outside a real codimension 2 locus $\Gamma\subset B_\infty$ (resp. $\check{\Gamma} \subset \check{B}_\infty$), $d_\infty$ (resp. $\check{d}_\infty$) is induced by a Monge-Amp\`ere metric; and
\item[(iv)]
the Monge-Amp\`ere manifolds $B_\infty\setminus \Gamma$ and $\check{B}_\infty\setminus \check{\Gamma}$ are Legendre dual to each other.
\end{itemize}

This limiting version of the SYZ conjecture was proposed independently by Gross-Wilson \cite{Gross-Wilson00} and Kontsevich-Soibelman \cite{Kontsevich-Soibelman01}. In fact, the general question of understanding the limiting behavior of Ricci-flat metrics was raised by Yau in his famous lists of open problems \cite{Yau82, Yau93}. Motivated by the SYZ picture of mirror symmetry, this question has been studied extensively in the last 15 years, and substantial progress has been made by Gross-Wilson \cite{Gross-Wilson00}, Tosatti \cite{Tosatti09, Tosatti10}, Ruan-Zhang \cite{Ruan-Zhang11}, Zhang \cite{Zhang09}, Rong-Zhang \cite{Rong-Zhang11, Rong-Zhang12} and more recently, Gross-Tosatti-Zhang \cite{GTZ13, GTZ13a}.

The metric spaces $B_\infty$ and $\check{B}_\infty$ should be thought of as limits of bases of SYZ fibrations on the families of Calabi-Yau manifolds. Applying the above picture, one may try to construct the mirror of a maximally unipotent degeneration of Calabi-Yau manifolds $\mathfrak{X} \to D$ as follows. We first identify the Gromov-Hausdorff limit $B_\infty$. Then we take the Legendre dual $\check{B}_0$ of $B_\infty \setminus \Gamma$ and try to compactify the quotient $\check{X}_0 := T\check{B}_0 / \Lambda$ to get the correct mirror. Unfortunately this na\"ive approach will not work because the natural complex structure on $\check{X}_0$ is not globally defined due to nontrivial monodromy of the affine structure around the singularities in $\check{B}_\infty$. To get the corrected mirror, one needs to deform the complex structure on $\check{X}_0$ by taking into account contributions from holomorphic disks. 

\section{SYZ for compact Calabi-Yau manifolds}



It is expected that the symplectic structure that we get using semi-flat SYZ mirror symmetry can naturally be compactified to give a global symplectic structure on the mirror. Indeed the work of Casta\~no-Bernard and Matessi \cite{Castano-Bernard-Matessi09} have shown that the topological Calabi-Yau compactifications constructed by Gross in \cite{Gross-Inventiones} can be made into symplectic compactifications, hence producing pairs of compact symplectic 6-folds which are homeomorphic to known mirror pairs of Calabi-Yau 3-folds (such as the quintic 3-fold and its mirror) and equipped with Lagrangian torus fibrations whose bases are Legendre dual integral affine manifolds with singularities.

On the other hand, as we have mentioned before, it was already anticipated in the original SYZ proposal \cite{SYZ96} that the Ricci-flat metric on the mirror should differ from the semi-flat Calabi-Yau metric \cite{GSVY90} by instanton corrections from holomorphic disks whose boundaries wrap non-trivial 1-cycles in the fibers of an SYZ fibration. Since the metric on the mirror is determined uniquely by its symplectic and complex structures, it is natural to expect that the instanton corrections are all contributing to perturbations of the complex structure on the mirror. This is indeed the key idea underlying the SYZ conjecture. As holomorphic disks can be glued to give holomorphic curves, this explains why mirror symmetry can lead to enumerative predictions.

Now, given an affine manifold with singularities $B$. Let $\Gamma \subset B$ be the singular locus and denote by $B_0 = B\setminus \Gamma$ the smooth part. Then one would like to construct a complex manifold which is a compactification of a small deformation of $X_0 := TB_0 / \Lambda$. This is called the {\em reconstruction problem}, which lies at the heart of the algebraic-geometric SYZ program of Gross and Siebert \cite{Gross-Siebert03, Gross-Siebert-logI, Gross-Siebert-logII, Gross-Siebert-reconstruction}.

The reconstruction problem was first attacked by Fukaya \cite{Fukaya05}. He attempted to find suitable perturbations by directly solving the Maurer-Cartan equation which governs the deformations of complex structures on $X_0$. In the dimension two case, his heuristic arguments suggested that the perturbations should come from gradient flow trees in $B$ whose ends emanate from the singular set $\Gamma$. The latter should be limits of holomorphic disks bounding Lagrangian fibers and singularities of an SYZ fibration when one approaches a large complex structure limit. Fukaya made a series of beautiful conjectures explaining how quantum corrections are contributing to the complex structure on the mirror and provided an intuitively clear picture, but unfortunately the analysis required to make his arguments rigorous seemed out of reach.

Kontsevich and Soibelman \cite{Kontsevich-Soibelman06} got around the analytic difficulties in Fukaya's arguments by working with rigid analytic spaces. They started with an integral affine structure on $S^2$ with 24 singular points such that the monodromy of the affine structure around each singular point is the simplest one:
$\left(\begin{array}{cc}
1 & 1\\
0 & 1
\end{array}\right)$, and managed to construct a non-Archimedean analytic $K3$ surface. The basic idea is to attach an automorphism to each gradient flow line in Fukaya's construction, and modify the gluing between complex charts in the mirror by these automorphisms, thereby resolving the incompatibilities between charts which arise from the nontrivial monodromy of the affine structure around the discriminant locus. A crucial step in their argument is a key lemma showing that when two gradient flow lines intersect, one can always add new lines together with new automorphisms attached so that the composition around each intersection point is the identity. This is called the {\em scattering phenomenon}, which assures that the composition of automorphisms attached to lines crossed by a path is independent of the path chosen, hence guaranteeing that the modified gluings are consistent.


At around the same time, Gross and Siebert launched their spectacular program \cite{Gross-Siebert03, Gross-Siebert-logI, Gross-Siebert-logII, Gross-Siebert-reconstruction} aiming at an algebraic-geometric approach to the SYZ conjecture. Motivated by the limiting version of the SYZ conjecture we discussed in the previous section and the observation by Kontsevich that the Gromov-Hausdorff limit will be roughly the dual intersection complex of the degeneration, they formulated an {\em algebraic-geometric SYZ procedure} to construct the mirror. In more details, starting with a toric degeneration of Calabi-Yau manifolds, the first step is to construct the dual intersection complex. Then one takes the (discrete) Legendre transform and try to reconstruct the mirror toric degeneration of Calabi-Yau manifolds from the Legendre dual. In this way, they can completely forget about SYZ fibrations. The claim is that all information are encoded in the tropical geometry of the dual intersection complex, which is an integral affine manifold with singularities and plays the role of the base of an SYZ fibration.

Using the above key lemma of Kontsevich and Soibelman, together with many new ideas such as employing log structures and techniques from tropical geometry, Gross and Siebert eventually succeeded in giving a solution to the reconstruction problem in {\em any} dimension \cite{Gross-Siebert-reconstruction}. More precisely, given any integral affine manifold with singularities satisfying certain assumptions and equipped with some additional structures like a polyhedral decomposition, they constructed a toric degeneration of Calabi-Yau manifolds which can be described explicitly and canonically via tropical trees in $B$. Furthermore, the Calabi-Yau manifolds they constructed are defined over $\bC$, instead of being rigid analytic spaces.

Now the goal is to acquire a conceptual understanding of mirror symmetry by going through the tropical world. On the B-side, Gross and Siebert conjectured that the deformation parameter in their construction is a canonical coordinate and period integrals of the family of Calabi-Yau manifolds can be expressed in terms of tropical disks in $B$. They have already given some evidences in the local cases (such as the local $\bP^2$ example in \cite[Remark 5.1]{Gross-Siebert-reconstruction}) and are working out the general case.

On the A-side, one would like to understand the Gromov-Witten theory of a smooth fiber by working entirely on the central fiber of a toric degeneration, whose dual intersection complex is precisely the affine manifold that Gross and Siebert started with. The recent independent works of Abramovich and Chen \cite{Chen10, Abravomich-Chen11} and Gross and Siebert \cite{Gross-Siebert13}, which developed the theory of {\em log Gromov-Witten invariants}, generalizing previous works of A.-M. Li and Ruan \cite{Li-Ruan01}, Ionel and Parker \cite{Ionel-Parker03, Ionel-Parker04}, and Jun Li \cite{Li02} on relative Gromov-Witten theory, constituted a significant step towards this goal. If furthermore one can prove a general correspondence theorem between tropical and holomorphic curves/disks, in the same vein as the works of Mikhalkin \cite{Mikhalkin05, Mikhalkin06}, Nishinou-Siebert \cite{Nishinou-Siebert06} and Nishinou \cite{Nishinou12, Nishinou09}, then we would be able to connect the A-side (i.e. Gromov-Witten theory) of a Calabi-Yau manifold to the tropical world.

Albeit much work needs to be done, this lays out a satisfying picture explaining the geometry of mirror symmetry via tropical geometry. We refer the reader to the beautiful survey articles of the inventors \cite{Gross-Siebert11a, Gross12} for overviews of the Gross-Siebert program.





\section{SYZ for noncompact Calabi-Yau manifolds}

The lack of examples of special Lagrangian torus fibrations is one main obstacle in implementing the original SYZ proposal for compact Calabi-Yau manifolds (and perhaps this is one of the reasons why Gross and Siebert wanted to develop an algebraic-geometric version). But there are plenty of noncompact examples where one can find explicit special Lagrangian torus fibrations, such as those constructed by Goldstein \cite{Goldstein01} and Gross \cite{Gross01} in the case of toric Calabi-Yau manifolds. Moreover, open Gromov-Witten invariants which count maps from open Riemann surfaces to the manifold are well-defined in the toric case by the works of Fukaya, Oh, Ohta and Ono \cite{FOOO-toricI, FOOO-toricII, FOOO-toricIII} (and even in the $S^1$-equivariant case by as has been done in the thesis of Liu \cite{Liu02}). So it makes perfect sense to carry out the SYZ proposal directly for toric Calabi-Yau manifolds, without retreating to the tropical world.

This brings us to the realm of {\em local mirror symmetry}, which was originally an application of mirror symmetry techniques to Fano surfaces within compact Calabi-Yau 3-folds, and could be derived using physical arguments from mirror symmetry for compact Calabi-Yau manifolds by taking certain limits in the K\"ahler and complex moduli spaces \cite{KKV97}. Since this mirror symmetry provides numerous interesting examples and predictions, it has been attracting much attention from both physicists and mathematicians \cite{Leung-Vafa98, CKYZ99, HIV00, Gross01, Gross-Inventiones, Takahashi01, Klemm-Zaslow01, Graber-Zaslow02, Hosono00, Hosono06, Forbes-Jinzenji05, Forbes-Jinzenji06, Konishi-Minabe10, Seidel10}.

Let $X$ be an $n$-dimensional toric Calabi-Yau manifold, which is necessarily noncompact. To carry out the SYZ construction, we consider a special Lagrangian torus fibration $\mu: X \to B$ constructed by Goldstein and Gross; such a fibration is {\em non-toric}, meaning that it is not the usual moment map associated to the Hamiltonian $T^n$-action on $X$. The discriminant locus of this SYZ fibration has been analyzed in details by Gross and can be described explicitly. Topologically, the base $B$ is simply an upper half-space in $\mathbb{R}^n$, and it admits an integral affine manifold with both singularities and boundary. The pre-image of the boundary $\partial B \subset B$ is a non-toric smooth hypersurface $D \subset X$. The discriminant locus $\Gamma \subset B$ is a real codimension two tropical subvariety sitting inside a hyperplane $H$ which we call the {\em wall} in $B$. By definition, the wall(s) inside the base of an SYZ fibration is the loci of Lagrangian torus fibers which bound holomorphic disks with Maslov index zero in $X$. It divides the base into different chambers over which the Lagrangian torus fibers behave differently in a Floer-theoretic sense. In the case of the Gross fibration, the wall $H \subset B$, which is parallel to the boundary hyperplane $\partial B$, divides the base into two chambers.

Now one considers (virtual) counts of holomorphic disks in $X$ bounded by fibers of the SYZ fibration $\mu$ which intersect with the hypersurface $D$ at one point with multiplicity one (i.e. disks with Maslov index two). As a point moves from one chamber to another across the wall, the virtual number of holomophic disks bounded by the corresponding Lagrangian fiber (or {\em genus 0 open Gromov-Witten invariants}) jumps, exhibiting a {\em wall-crossing phenomenon}. This has been analyzed by Auroux \cite{Auroux07, Auroux09} and Chan, Lau and Leung \cite{CLL12} by applying the sophisticated machinery developed by Fukaya, Oh, Ohta and Ono \cite{FOOO-book}.  Notice that there is no scattering phenomenon in this case because there is only one wall. By applying the SYZ dual fibration construction on each chamber in the base, and then gluing the resulting pieces together according to the wall-crossing formulas, we obtain the {\em instanton-corrected} or {\em SYZ mirror family} $\check{X}$, which is parametrized by the K\"ahler moduli space of $X$ \cite{CLL12, AAK12}. The result agrees with the predictions by physical arguments \cite{CKYZ99, HIV00}.

This SYZ mirror construction is very precise in the sense that it tells us exactly which complex structure on $\check{X}$ is corresponding to any given symplectic structure on $X$ -- the defining equation of the mirror $\check{X}$ is an explicit expression written entirely in terms of the K\"ahler parameters and disk counting invariants of $X$. For example, the SYZ mirror of $X=K_{\bP^2}$ is given by\footnote{More precisely, the SYZ mirror of $K_{\bP^2}$ is the Landau-Ginzburg model $(\check{X}, W=u)$; see the next section.}
\begin{equation}\label{eqn:mirror_KP2}
\check{X} = \left\{(u,v,z_1,z_2)\in\bC^2\times(\bC^\times)^2 \mid uv = 1+\delta(q) + z_1 + z_2 + \frac{q}{z_1z_2} \right\},
\end{equation}
where $q$ is the K\"ahler parameter measuring the symplectic area of a projective line inside the zero section of $K_{\bP^2}$ over $\bP^2$, and
\begin{equation}\label{eqn:gen_fcn_openGW_KP2}
1+\delta(q) = \sum_{k=0}^\infty n_k q^k
\end{equation}
is a generating series of genus 0 open Gromov-Witten invariants.

Furthermore, the SYZ construction naturally defines the {\em SYZ map}, which is a map from the K\"ahler moduli space of $X$ to the complex moduli space of $\check{X}$. As conjectured by Gross and Siebert, the SYZ mirror family should be written in canonical coordinates, or put it in another way, the SYZ map should give an inverse to the mirror map. Evidences for this conjecture for toric Calabi-Yau surfaces and 3-folds were given in \cite{CLL12, LLW12}, and Chan, Lau and Tseng \cite{CLT11} proved the conjecture in the case when $X$ is the total space of the canonical line bundle over a compact toric Fano manifold. Recently, by applying orbifold techniques, the conjecture was proved for {\em all} toric Calabi-Yau manifolds in \cite{CCLT13}.

The main challenge in proving these results is the computation of the genus 0 open Gromov-Witten invariants defined by Fukaya, Oh, Ohta and Ono \cite{FOOO-toricI}. Since the moduli spaces of holomorphic disks are usually highly obstructed, these invariants are in general very difficult to compute. Currently, there are only very few techniques available (such as open/closed equalities, toric mirror theorems, degeneration techniques, etc). For example, the invariants in \eqref{eqn:gen_fcn_openGW_KP2} can be computed:
\begin{align*}
n_k = 1, -2, 5, -32, 286, -3038, 35870, \ldots
\end{align*}
for $k=0,1,2,3,4,5,6,\ldots$, which agrees with period computations in \cite{Graber-Zaslow02}.

We should mention that the SYZ construction can be carried out also in the reverse direction \cite{AAK12} (see also \cite[Section 5]{Chan12}). For example, starting with the conic bundle \eqref{eqn:mirror_KP2}, one can construct an SYZ fibration using similar techniques as in \cite{Goldstein01, Gross01}. Although the discriminant locus is of real codimension one in the case, one can construct the SYZ mirror and this gives us back the toric Calabi-Yau 3-fold $K_{\bP^2}$, as expected.\footnote{More precisely, the SYZ mirror of \eqref{eqn:mirror_KP2} is the complement of a smooth hypersurface in $K_{\bP^2}$.}

Nevertheless, outside the toric setting, it is not clear how SYZ constructions can be performed in such an explicit way. One major problem is the well-definedness of open Gromov-Witten invariants. Only in a couple of non-toric cases (see Liu \cite{Liu02} and Solomon \cite{Solomon06}) do we have a well-defined theory of open Gromov-Witten invariants.\footnote{There are, however, recent works of Fukaya \cite{Fukaya10, Fukaya11} on defining disk counting invariants for compact Calabi-Yau 3-folds.}





\section{SYZ in the non--Calabi-Yau setting}

Not long after its discovery, mirror symmetry has been extended to the non--Calabi-Yau setting, notably to Fano manifolds, through the works of Batyrev \cite{Batyrev93}, Givental \cite{Givental95, Givental96, Givental98}, Kontsevich \cite{Kontsevich-ENS98}, Hori-Vafa \cite{Hori-Vafa00} and many others. Unlike the Calabi-Yau case, the mirror is no longer given by a manifold; instead, it is predicted to be a pair $(\check{X},W)$, where $\check{X}$ is a non-compact K\"ahler manifold and $W:\check{X} \to \bC$ is a holomorphic function. In the physics literature, such a pair $(\check{X},W)$ is called a {\em Landau-Ginzburg model}, and $W$ is called the {\em superpotential} of the model \cite{Vafa91, Witten93}.

It is natural to ask whether the SYZ proposal continues to work in this setting as well. Auroux \cite{Auroux07} was the first to consider this question and in fact he extended the SYZ proposal to a much more general setting. Namely, he considered pairs $(X,D)$ consisting of a compact K\"ahler manifold $X$ together with an effective anticanonical divisor $D$. The defining section of $D$ gives a holomorphic volume form on $X\setminus D$ with simple poles along the divisor $D$, so it makes sense to speak about special Lagrangian torus fibrations on the complement $X\setminus D$. Suppose that we are given such a fibration $\mu:X\setminus D \to B$, then we can try to produce the SYZ mirror $\check{X}$ by $T$-duality (i.e. consider the moduli space of pairs $(L,\nabla)$ where $L$ is a fiber of $\mu$ and $\nabla$ is a flat $U(1)$-connection over $L$) modified by instanton corrections. Moreover, the superpotential $W$ will naturally appears as the object mirror to Fukaya-Oh-Ohta-Ono's obstruction chain $\mathfrak{m}_0$.

When $X$ is a compact toric K\"ahler manifold, a canonical choice of $D$ is the union of all toric prime divisors. Also, the moment map provides a convenient Lagrangian torus fibration on $X$, which has the nice property that it restricts to a torus bundle on the open dense torus orbit $X\setminus D$. In this case, the SYZ mirror manifold $\check{X}$ is simply given by the algebraic torus $(\bC^\times)^n$, because we have a torus bundle and there are no instanton corrections in the construction of the mirror manifold. All the essential information is encoded in the superpotential $W$. Prior to the work of Auroux, it was Cho and Oh \cite{Cho04, Cho-Oh06} who first noticed that $W$ can be expressed in terms disk counting invariants (or open Gromov-Witten invariants). By classifying all holomorphic disks in $X$ bounded by moment map fibers, they got an explicit formula for $W$ in the case when $X$ is Fano, and this agrees with earlier predictions obtained using physical arguments by Hori and Vafa \cite{Hori-Vafa00}. This was later vastly generalized by the works of Fukaya, Oh, Ohta and Ono \cite{FOOO-toricI, FOOO-toricII, FOOO-toricIII} on Lagrangian Floer theory and mirror symmetry for toric manifolds.

In \cite{Chan-Leung10a}, mirror symmetry for toric Fano manifolds was used as a testing ground to see how useful Fourier-Mukai--type transforms, or what we call SYZ transforms, could be in the investigation of the geometry of mirror symmetry. For a toric Fano manifold $X$, we consider the open dense torus orbit $X_0 := X\setminus D \subset X$, which is also the union of Lagrangian torus fibers of the moment map. Symplectically, we can write $X_0 = T^*B_0/\Lambda^\vee $, where $B$ is the moment polytope and $B_0$ denotes its interior. Then the SYZ mirror is $\check{X} := TB_0/\Lambda$ which is a bounded domain in $(\bC^\times)^n$. To obtain the superpotential $W$, we consider the space
$$\tilde{X} := X_0 \times \Lambda \subset \mathcal{L} X$$
of fiberwise geodesic/affine loops in $X$. On $\tilde{X}$, we have an instanton-corrected symplectic structure $\tilde{\omega} = \omega + \Phi$, where $\Phi$ is a generating function of genus 0 open Gromov-Witten invariants which count (virtually) holomorphic disks bounded by moment map fibers.

An explicit SYZ transform $\mathcal{F}$ was then constructed by combining the semi-flat SYZ transform $\mathcal{F}^\textrm{semi-flat}$ with fiberwise Fourier series, and it was shown that $\mathcal{F}$ transforms the corrected symplectic structure $\tilde{\omega}$ on $X$ precisely to the holomorphic volume form $e^W\check{\Omega}$ of the mirror Landau-Ginzburg model $(\check{X},W)$, where $W$ was obtained by taking fiberwise Fourier transform of $\Phi$. Moreover, $\mathcal{F}$ induces an isomorphism between the (small) quantum cohomology ring $QH^*(X)$ of $X$ and the Jacobian ring $Jac(W)$ of $W$. The proof was by passing to the tropical limit, and observing that a tropical curve whose holomorphic counterpart contributes to the quantum product can be obtained as a gluing of tropical disks (see the work \cite{Chan-Leung10a} for more details). This observation was later generalized and used by Gross \cite{Gross10} in his study of mirror symmetry for the big quantum cohomology of $\bP^2$ via tropical geometry.


As for manifolds of general type, there are currently two main approaches to their mirror symmetry along the SYZ perspective. One is the work of Abouzaid, Auroux and Katzarkov, where they considered a hypersurface $H$ in a toric variety $V$ and constructed a Landau-Ginzburg model which is SYZ mirror to the blowup of $V\times \bC$ along $H\times \{0\}$. In particular, when $H$ is the zero set of a bidegree $(3,2)$ polynomial in $V = \bP^1\times\bP^1$, their construction produces a mirror of the genus 2 Riemann surface, which is in agreement with a previous proposal by Katzarkov \cite{Katzarkov07, KKOY09, Seidel11}.

Another approach, which is more in line with the Gross-Siebert program, is the work by Gross, Katzarkov and Ruddat \cite{GKR12}. They proposed that the mirror to a variety of general type is a reducible variety equipped with a certain sheaf of vanishing cycles. Presumably, the mirror produced in this approach should give the same data as the one produced by \cite{AAK12}. For example, the reducible variety should the critical locus of the superpotential of the SYZ mirror Landau-Ginzburg model. But the precise relations between these two approaches are still under investigation.








\section{Beyond SYZ}

Besides providing a beautiful geometric explanation of mirror symmetry, the SYZ conjecture \cite{SYZ96} has been exerting its long-lasting effect on many related areas of mathematics as well. Let us briefly describe several examples of applications in this regard.\\

\noindent\textbf{HMS via SYZ.} As we have seen, the SYZ conjecture is based upon the idea of D-branes in string theory. Recall that B-branes (i.e. D-branes in the B-model) are coherent sheaves over complex subvarieties while A-branes (i.e. D-branes in the A-model) are special Lagrangian submanifolds equipped with flat $U(1)$ connections. It therefore makes sense to view Kontsevich's HMS conjecture \cite{Kontsevich-ICM94}, which asserts that the Fukaya category of a Calabi-Yau manifold $X$ is equivalent to the derived category of coherent sheaves on the mirror $\check{X}$, as a manifestation of the isomorphism between the A-model on $X$ and the B-model on $\check{X}$. So rather naturally, one expects that the SYZ proposal, and in particular SYZ transforms, can be exploited to construct functors which realize the categorial equivalences asserted by the HMS conjecture.

For example, given a Lagrangian section of a special Lagrangian torus fibration $\mu:X \to B$, its intersection point with a fiber $L$ of $\mu$ determines a flat $U(1)$-connection on the dual torus $L^\vee$. Patching these flat $U(1)$-connections together should give a holomorphic line bundle over the total space of the dual fibration, which is the mirror $\check{X}$. This simple idea, first envisioned by Gross \cite{Gross98, Gross99}, was explored by Arinkin and Polishchuk \cite{Arinkin-Polishchuk01} and Leung, Yau and Zaslow \cite{LYZ00} to construct SYZ transforms, which were applied to prove and understand the HMS conjecture in the semi-flat Calabi-Yau case. Later, the same idea was also employed to study the HMS conjecture for toric varieties \cite{Abouzaid06, Abouzaid09, Fang08, FLTZ12, FLTZ11b, Chan09, Chan-Leung12, CHL12}.

In some more recent works \cite{Chan12, Chan-Ueda12, CPU13}, SYZ transforms were applied to construct geometric Fourier--type functors (on the objects level) which realize the HMS categorial equivalences for certain examples of local Calabi-Yau such as resolutions of the $A_n$-singularities and the smoothed conifold, where one encounters SYZ fibrations with singular fibers and hence nontrivial quantum corrections. On the other hand, work in progress by K.-L. Chan, Leung and Ma \cite{CLM12, CLM13} have shown that SYZ transforms can also used to construct the HMS equivalences on the morphism level, at least in the semi-flat case. The ultimate goal is to construct a canonical geometric Fourier-type functor associated to any given SYZ fibration, which realizes the equivalences of categories asserted by the HMS conjecture, thereby enriching our understanding of the geometry of the HMS conjecture, and also mirror symmetry as a whole.\\

\noindent\textbf{Ricci-flat metrics and disk counting.} A remarkable observation in the SYZ paper \cite{SYZ96} is that a Ricci-flat metric on the mirror can be decomposed as the sum of a semi-flat part (which was written down explicitly earlier in \cite{GSVY90}) and an instanton-corrected part which should come from contributions by holomorphic disks in the original Calabi-Yau manifold which have boundaries on the Lagrangian torus fibers of an SYZ fibration. This suggests a concrete and qualitative description of Ricci-flat metrics, which are extremely hard to write down.

In general, such a qualitative description is still highly nontrivial to obtain because it is difficult to find nontrivial examples of SYZ fibrations on compact Calabi-Yau manifolds, and open Gromov-Witten theory is not yet well-understood. However, recent pioneering works of Gaiotto, Moore and Neitzke \cite{GMN10, GMN13} have shed new light on the hyperk\"ahler case. They proposed a new (partially conjectural) relation between hyperk\"ahler metrics on the total spaces of complex integrable systems (the simplest example of which reproduces the well-known Ooguri-Vafa metric \cite{Ooguri-Vafa96}) and Kontsevich-Soibelman's wall-crossing formulas.

To describe their proposal in a bit more details, let us consider a complex integrable system $\psi: M \to B$, i.e. $M$ is holomorphic symplectic and the fibers of $\psi$ are complex Lagrangian submanifolds. (More precisely, what Gaiotto, Moore and Neitzke were looking at in \cite{GMN10, GMN13} were all meromorphic Hitchin systems, in which case complete hyperk\"ahler metrics were first constructed by Biquard and Boalch \cite{Biquard-Boalch04}; see \cite[Section 4.1]{GMN13}.) They made use of the fact that any hyperk\"ahler metric is characterized by the associated twistor space, and tried to write down a $\bC^\times$-family of holomorphic Darboux coordinates on $M$ which satisfy the hypotheses of the theorem of Hitchin et al. \cite{HKLR87}. In particular, they required the coordinates to satisfy certain wall-crossing formulas which describe the discontinuity of the coordinates across the so-called {\em BPS rays}, where the (virtual) counts of BPS states jump.

These wall-crossing formulas turn out to be of the same kind as those used by Kontsevich and Soibelman \cite{Kontsevich-Soibelman06} and Gross and Siebert \cite{Gross-Siebert-reconstruction} in their constructions of toric degenerations of Calabi-Yau manifolds (and on the other hand they are the same as the wall-crossing formulas in motivic Donaldson-Thomas theory \cite{Joyce-Song12, Kontsevich-Soibelman08}). In view of this and the SYZ conjecture (and also the fact that a vast family of examples of noncompact SYZ fibrations on meromorphic Hitchin systems, including many in complex dimension two (e.g. gravitational instantons, log-Calabi-Yau surfaces) have been constructed via hyperk\"ahler rotation in \cite{Biquard-Boalch04}), it is natural to expect that the hyperk\"ahler metrics on those complex integrable systems considered by Gaiotto, Moore and Neitzke can be expressed in terms of holomorphic disks.

This was done for the simplest example -- the Ooguri-Vafa metric in \cite{Chan10}. More recent works of W. Lu \cite{Lu10, Lu12} have demonstrated that in general the twistor spaces and holomorphic Darboux coordinates on meromorphic Hitchin systems studied in \cite{GMN10, GMN13} produced the same data as those required to run the Gross-Siebert program \cite{Gross-Siebert-reconstruction}, hence showing that there must be some (perhaps implicit) relations between the hyperk\"ahler metrics and tropical disks counting. In his PhD thesis \cite{Lin-PhDthesis}, Y.-S. Lin considered elliptic $K3$ surfaces and tropical disk counting invariants. He proved that his invariants satisfy the same wall-crossing formulas as appeared in \cite{GMN10, GMN13}. This again shows that the hyperk\"ahler metrics on those K3 surfaces are closely related to disk counting. There are also recent works by Stoppa and his collaborators \cite{Stoppa11, Filippini-Stoppa13} demonstrating the intimate relations between the wall-crossing formulas in motivic Donaldson-Thomas theory and the constructions of Gaiotto, Moore and Neitzke.\\

\noindent\textbf{Other applications of SYZ.} Let us also mention two recent, unexpected applications of SYZ constructions.

In their recent joint project \cite{GHK11}, Gross, Hacking and Keel constructed mirror families to log Calabi-Yau surfaces, i.e. pairs $(Y,D)$ where $Y$ is a nonsingular projective rational surface and $D \in |-K_Y|$ is a cycle of rational curves, by extending the construction in \cite{Gross-Siebert-reconstruction} to allow the affine manifolds to have more general (i.e. worse) singularity types. Amazingly, their results could be applied to give a proof of a 30-year-old conjecture of Looijenga \cite{Looijenga81} concerning smoothability of cusp singularities.

In an even more recent preprint \cite{GHK12}, they applied their construction again to prove a Torelli theorem for log Calabi-Yau surfaces, which was originally conjectured in 1984 by Friedman \cite{Friedman12}. On the other hand, their construction is also closely connected with the theory of cluster varieties, and they have suggested a vast generalisation of the Fock-Goncharov dual bases. For a nice exposition of these exciting new results and developments, we refer the reader to the nice survey article by Gross and Siebert \cite{Gross-Siebert12}.

In another unexpected direction, the SYZ construction has recently been applied to construct new {\em knot invariants}. For a knot $K$ in $S^3$, its conormal bundle $N^*K$ is canonically a Lagrangian cycle in the cotangent bundle $T^*S^3$. In \cite{DSV11}, Diaconescu, Shende and Vafa constructed a corresponding Lagrangian cycle $L_K$ in the resolved conifold $X := \mathcal{O}_{\bP^1}(-1)\oplus\mathcal{O}_{\bP^1}(-1)$, which is roughly speaking done by lifting the conormal bundle $N^*K$ off the zero section and letting $T^*S^3$ undergo the conifold transition. Their construction was motivated by a mysterious phenomenon called {\em large $N$ duality} in physics.

In \cite{Aganagic-Vafa12}, Aganagic and Vafa defined a new knot invariant by a generalized SYZ construction applied to the pair $(X,L_K)$. Roughly speaking, their invariant is the generating series of open Gromov-Witten invariants for $(X,L_K)$. It turned out that the resulting invariant is always a polynomial and they conjectured that it should be a deformation of the classical A-polynomial in knot theory \cite{CCGLS94}. Furthermore, an interesting relation between their invariant and augmentations of the contact homology algebra of $K$ \cite{Ng12} was suggested. Substantial evidences for this relation was obtained in a very recent preprint \cite{AENV13}.

These two new applications of the SYZ conjecture, together with many more which are yet to come, open up new directions in mirror symmetry and other branches of mathematics and physics,\footnote{One interesting story that we have not mentioned is the application of the SYZ picture to $G_2$ manifolds by Gukov, Yau and Zaslow \cite{GYZ03} which was aimed at explaining the duality between M-theory and heterotic string theory.} and they are all pointing out to further research works for the future.

\bibliographystyle{amsplain}
\bibliography{geometry}

\end{document}